\documentclass{article}
\usepackage{latexsym}
\begin{document}

\title{Boundary value problems in dimensions seven, four and three related to exceptional holonomy}
\author{Simon Donaldson}
\maketitle


\newcommand{\bC}{{\bf C}}
\newcommand{\bP}{{\bf P}}
\newcommand{\bR}{{\bf R}}
\newcommand{\Bix}{\Box}
\newcommand{\trho}{\hat{\rho}}
\newcommand{\urho}{\underline{\rho}}
\newcommand{\uomega}{\underline{\omega}}
\newcommand{\umu}{\underline{\mu}}
\newcommand{\wumu}{\hat{\underline{\mu}}}
\newcommand{\uh}{\underline{h}}
\newtheorem{lem}{Lemma}
\newtheorem{prop}{Proposition}
\newcommand{\hook}{\leftharpoonup}
\ \ \ \ \ \ \ \ \ \ \ \ \ \ \ \ {\em Dedicated to Nigel Hitchin, for his 70th. birthday.}

\

The variational point of view on exceptional structures in dimensions 6,7 and 8 is one of Nigel Hitchin's seminal contributions. One feature of this point of view is that it motivates the study of boundary value problems, for structures with prescribed data on a boundary.  In this article we consider the case of 7 dimensions and \lq\lq $G_{2}$-manifolds''. We will review briefly a general framework and then go on to examine in more detail symmetry reductions to dimensions 4 (in Section 2) and 3 (in Section 3). In the latter case we encounter an interesting variational problem related to the real Monge-Amp\`ere equation and in Section 4 we describe a generalisation of this.

 The author is grateful to Claude Le Brun and Lorenzo Foscolo for valuable comments and discussions. 

 \section{The volume functional in 7 dimensions}
 
 Let $V$ be a 7-dimensional oriented real vector space. A $3$-form $\phi\in \Lambda^{3}V^{*}$ defines a quadratic form on $V$ with values in the real line $\Lambda^{7}V^{*}$ by the assignment
$$  v\mapsto (i_{v}\phi)^{2} \wedge \phi. $$
The fixed orientation means that it makes sense to say that this form is positive definite, and in that case we call $\phi$ a {\it positive} $3$-form. From the definition, a positive $3$-form defines a conformal class of Euclidean structures on $V$ and the ambiguity of scale can be fixed by choosing the Euclidean structure so that $\vert \phi\vert^{2}=7$, using the standard induced metric on $\Lambda^{3}V^{*}$. 

Now let $M$ be an oriented 7-manifold. Applying the above in each tangent space, we have the notion of a positive 3-form $\phi\in \Omega^{3}(M)$ and such a form defines a metric $g_{\phi}$ and volume form $\nu_{\phi}\in \Omega^{7}(M)$.
The first variation of the volume form with respect to a variation $\delta\phi$ in $\phi$ is given by
$$ \delta \nu_{\phi}= \delta \phi \wedge \Psi  $$
where $\Psi=\Psi(\phi)$ is a $4$-form determined by $\phi$,  which can also be expressed as
$$ 3 \Psi(\phi)= *_{g_{\phi}} \phi. $$
 Suppose that $M$ is a closed $7$-manifold and that $c\in H^{3}(M,\bR)$ is a cohomology class which can be represented by  positive $3$-forms, so we have a non-empty set ${\cal S}_{c}\subset \Omega^{3}(M)$ of  closed positive forms representing $c$. Hitchin's idea (\cite{kn:H1}, \cite{kn:H2}) is to consider the total volume
\begin{equation}   {\rm Vol}(\phi)= \int_{M} \nu_{\phi} \end{equation}
as a functional on ${\cal S}_{c}$. The first variation, with respect to a variation $\delta \phi=da$,  in $\phi$ is
\begin{equation}   \int_{M} \Psi\wedge da = - \int_{M} d\Psi \wedge a \end{equation}
so the Euler-Lagrange equation defining critical points is
$$   d\Psi=0. $$
By a well-known result of Fern\'andez and Gray, the two equations $d\phi=0, d*_{g_{\phi}}\phi=0$ imply that the $G_{2}$ structure defined by $\phi$ is torsion-free, or equivalently that the metric $g_{\phi}$ has holonomy contained in $G_{2}$. So, from this point of view, the search for these special structures can be divided into two stages:
\begin{itemize}\item Identify manifolds $M$ and classes $c\in H^{3}(M)$ such that ${\cal S}_{c}$ is non-empty;
\item Study the variational problem for the volume functional on ${\cal S}_{c}$.
\end{itemize}

The local theory of such critical points, with respect to small variations in $c$ and $\phi$, is well-understood. Hitchin proved that any critical point is a {\it local maximum} and in fact a strict local maximum modulo diffeomorphisms. The proof is an application of Hodge theory. One of the many interesting and fundamental questions in this area is whether it is a {\it global maximum} over the whole space ${\cal S}_{c}$. Another standard fact (proved earlier by Bryant),  is that critical points are stable with respect to variations in the cohomology class $c$:  a critical point $\phi$ belonging to a class $c$ can be deformed to a critical point for nearby classes in $H^{3}(M)$. That is, the moduli space of $G_{2}$-structures is locally modelled on $H^{3}(M)$.

Now we introduce our boundary value problem. Let $M$ be a compact oriented $7$-manifold with boundary an oriented $6$-manifold $N=\partial M$. There is a similar notion of a positive $3$-form $\rho$ on $N$: this is just the condition that at each point $p\in N$ the form $\rho$ can be extended to a positive form on $TM_{p}$. A basic algebraic fact is that such a  positive 3-form in $6$-dimensions is equivalent to a reduction to $SL(3,\bC)$, that is, to  an almost-complex structure with a trivialisation of the \lq\lq canonical line bundle''. Fix a closed positive $3$-form $\rho$ on $N$. We  assume that the class $[\rho]\in H^{3}(N)$ is in the image of the restriction map from $ H^{3}(M)$. Define an {\it enhancement} of $\rho$ to be an equivalence class of closed $3$-forms on $M$ extending $\rho$,  under the equivalence relation $\phi_{1}\sim \phi_{1} +da$ where $a$ vanishes on $N$. Thus the set of enhancements of $\rho$ is an affine space modelled on $H^{3}(M,N)$. Fix an enhancement $\trho$ of $\rho$ and let ${\cal S}_{\trho}$ be the set of positive forms on $M$ in this  equivalence  class.
Suppose that ${\cal S}_{\trho}$ is nonempty and consider the volume functional on this set, just as before. For  a variation $\delta\phi=da$ with $a$ vanishing on the boundary the integration by parts (2) is still valid and the critical points are given by solutions of $d\Psi=0$ just as before. That is,  we are studying $G_{2}$-structures on $M$ with the given boundary value $\rho$ and in the given enhancement class. So we have the same two questions as before: identify enhanced boundary values $\trho$ such that ${\cal S}_{\trho}$ is non-empty and then study the variational problem.

We will not enter into a proper  discussion of the local theory of this boundary value problem (with respect to small variations in $\phi$ and $\trho$) here,  but we  make two simple observations. For the first, we say that a $G_{2}$-structure $\phi$ on a manifold $M$ with boundary $N$ is a {\it formal maximum} of the volume functional if for any $a\in \Omega^{2}(M)$ whose restriction to the boundary vanishes we have
      $$ \frac{d^{2}}{dt^{2}} {\rm Vol}(\phi+ t da) \leq 0$$
      at $t=0$. In other words, the Hessian of the volume functional is non-negative.
\begin{prop}
   Suppose that $M$ is the closure of a domain in a closed $G_{2}$-manifold $M^{+}$ with $3$-form $\phi$. Let $\trho$ be the enhanced boundary value given by the restriction of $\phi$ to $ N=\partial M $ and $M$.   Then $\phi\vert_{M}$ is a formal  maximum of the volume functional  on ${\cal S}_{\trho}$.
\end{prop}

Let $f$ be a defining function for $\partial M=f^{-1}(0)$, positive on the interior of $M$. Then any $2$-form $a$ on $M$ whose restriction to $\partial M$ is zero can be written as
$a= b+ \eta\wedge df$ where $b$ vanishes in $TM\vert_{\partial M}$. For small $\epsilon$,  let $\chi_{\epsilon}:M\rightarrow \bR $ be the composite of $f$ with a standard cut-off function, such that $\chi_{\epsilon}$ vanishes when $f\leq \epsilon$, is equal to $1$ when $f\geq \epsilon$ and with $\vert d \chi_{\epsilon}\vert \leq C \epsilon^{-1}$. Set $a_{\epsilon} = \chi_{\epsilon} a$. Then 
$$   d(a_{\epsilon}) = d \chi_{\epsilon}\ \wedge b$$
satisfies a uniform $L^{\infty}$ bound, independent of $\epsilon$ (since $b$ is $O(\epsilon)$ on the support of $d \chi_{\epsilon}$). It follows that
$$ \frac{d^{2}}{d t^{2}} {\rm Vol}(\phi+ t d a) = \lim_{\epsilon\rightarrow 0} \frac{d^{2}}{d t^{2}} {\rm Vol}(\phi+ t d a_{\epsilon}), $$
and the latter is non-positive since $a_{\epsilon}$ can be extended by zero over the closed manifold $M^{+}$ and then we can apply Hitchin's result.  

 In the other direction, critical points are not always {\it strict} local maxima, modulo diffeomorphisms. To give an example of this, we define for $v\in \bR^{7}$ with $\vert v\vert <1/2$ the manifold-with-boundary $M_{v}\subset \bR^{7}$ to be
$$   M_{v}= \overline{B}^{7}\setminus (v+ \frac{1}{2} B^{7}), $$ where $B^{7}$ is the open unit ball. Let $\phi_{0}$ be the standard flat $G_{2}$ structure on $\bR^{7}$ and let $\rho_{v}$ be its restriction to the boundary of $M_{v}$. In this case $H^{3}(M_{v}, \partial M_{v})=0$ so there is no extra enhancement data.
We can choose diffeomorphisms $F_{v}: M_{0}\rightarrow M_{v}$ such that
the restriction to the boundaries pulls back $\rho_{v}$ to $\rho_{0}$. Then $F_{v}^{*}(\phi_{0})$ are critical points for the boundary value problem on $M_{0}$ which are not all equivalent, by diffeomorphisms of $M_{0}$,  to $\phi_{0}$.

\section{Reduction to dimension 4.}

In this section we consider an interesting reduction of the $7$-dimensional theory to 4-dimensions, as follows. Take $M=X\times \bR^{3}$ where $X$ is an oriented $4$-manifold and consider $3$-forms of the shape
\begin{equation} \phi= \sum_{i=1}^{3} \omega^{i} d\theta_{i} - d\theta_{1}d\theta_{2}d\theta_{3} \end{equation}
where $\theta_{i}$ are co-ordinates on $\bR^{3}$ and $\omega^{i}$ are $2$-forms on $X$. The condition that $\phi$ is a positive $3$-form goes over to the condition that $(\omega^{i})$ form a \lq\lq positive triple", by which we mean that  at each point they span a maximal positive subspace for the wedge product form on $\Lambda^{2}T^{*}X$. More invariantly, we are considering positive forms $\phi$ which are preserved by the translation action of the $\bR^{3}$ factor and such that the orbits are \lq\lq associative'' submanifolds. The condition that $\phi$ be closed goes over to the condition that the $\omega^{i}$ are closed $2$-forms on $X$,  making up a \lq\lq hypersympletic'' structure.  These structures are of considerable interest in 4-dimensional geometry, see for example \cite{kn:D}, \cite{kn:FY}. 

Given such a triple $\uomega= (\omega^{i})$ we define  a volume form $\chi$ on $X$ by the following procedure. Let $\chi_{0}$ be any volume form and define a matrix $(q^{ij})$ by 
$$  \chi_{0} q^{ij}=   \omega^{i}\wedge \omega^{j}. $$
Now put
$$  \chi= \det(q^{ij})^{1/3} \chi_{0} . $$
It is clear that this is independent of the choice of $\chi_{0}$. The $7$-dimensional volume form associated to $\phi$ is $- \chi d\theta_{1} d\theta_{2}d\theta_{3}$. All our constructions will be invariant under the action of $SL(3,\bR)$ on $\bR^{3}$ so it will sometimes be clearer to introduce a $3$-dimensional oriented vector space $W$ with fixed volume element and consider our data $\uomega$ as an element  of $W\otimes \Omega^{2}(X)$. Then a choice of co-ordinate system on $W$ gives the description as a triple $(\omega^{1}, \omega^{2}, \omega^{3})$.

Given a positive triple $\omega^{i}$, we define a matrix $(\lambda^{ij})$ by
\begin{equation}    \omega^{i}\wedge \omega^{j} = \lambda^{ij} \chi. \end{equation}
Thus $\det(\lambda^{ij})=1$, by the definition of $\chi$. Write $(\lambda_{ij})$ for the inverse matrix and set
\begin{equation}  \Theta_{i} = \sum_{j=1}^{3} \lambda_{ij} \omega^{j}\end{equation}   The $4$-form defined by $\phi$ is
$$\Psi= \sum_{{\rm cyclic}} \Theta_{i} d\theta^{j}d\theta^{k}  + \chi, $$
where the notation means that $(ijk)$ runs over the three cyclic permutations of $(123)$. Thus the condition that a closed triple $(\omega^{i})$ defines a $G_{2}$ structure is $d\Theta_{i}=0$ which is to say:
\begin{equation}  \sum_{j=1}^{3} d\lambda_{ij} \wedge \omega^{j} =0.  \end{equation}
These equations are obviously satisfied when the matrix $(\lambda^{ij})$ is constant on $X$ and these solutions are the {\it hyperk\"ahler metrics}. 
Of course we can produce these equations (8) from a 4-dimensional reduction of Hitchin's variation formulation: the equations are the Euler-Lagrange equation for the functional
\begin{equation}  {\rm Vol}(\uomega)= \int_{X} \chi \end{equation}
on closed positive triples $\uomega$, with respect to exact variations of compact support.

 It is well-known, and easy to show directly, that the only solutions of the equations (6) on a compact 4-manifold are hyperk\"ahler and this gives extra motivation for considering the boundary value problem. So let $X$ be a $4$-manifold with boundary $Y$ and consider  triples $\umu=(\mu^{1},\mu^{2},\mu^{3})$ of closed $2$-forms on $Y$ which form a basis for $\Lambda^{2}T^{*}Y$ at each point. In our more invariant set-up, $\umu$ lies in $W\otimes \Omega^{2}(Y)$.  We define an enhancement $\umu$ in the obvious way, so the space of enhancements of a given $\umu$ is an affine space modelled on $W\otimes H^{2}(X,Y)$. Fix an enhancement $\wumu$ and let ${\cal S}_{\wumu}$ be the set of closed positive triples on $X$ in the given equivalence class. So the reduced versions of our questions are, first,  whether this  set is non-empty and, second, to study the variational problem given by the volume functional (7).

Stokes' Theorem implies that the integrals
$$   Q^{ij}= \int_{X} \omega^{i}\wedge \omega^{j}, $$
are independent of the choice of $\omega^{i}$ in a fixed enhancement class $\wumu$.
More invariantly, $Q$ is a quadratic form on our vector space $W$ and $\det Q$ is defined, as a real number, using the fixed volume form on $W$. 
This has two simple consequences.

\begin{prop}
  If ${\cal S}_{\wumu}$ is non-empty then $Q$ is positive definite and there is an upper bound
$$ \int_{X} \chi \leq  \det Q, $$
for $\uomega\in {\cal S}_{\wumu}$ and $\chi=\chi(\uomega)$. Equality holds if and only if $\uomega$ is hyperk\"ahler. 
\end{prop}
To see that $Q$ is positive definite it suffices, by change of basis, to see that $Q^{11}>0$. But this clear since $\omega^{1}\wedge\omega^{1}$ is positive pointwise on $X$. To establish the upper bound it suffices, by change of basis, to consider the case when $Q^{ij}=\delta^{ij}$. Recall that we write
$\omega^{i}\wedge \omega^{j}=\lambda^{ij} \chi$ where $\det(\lambda)=1$. Then we have the elementary inequality (the arithmetic-geometric mean inequality for the eigenvalues) ${\rm Tr}(\lambda)\geq 3$. So
$$ 3 \int_{X}\chi\leq  \int_{X}{\rm Tr} \lambda \ \chi  =\sum_{i} \int_{X}(\omega^{i})^{2} = 3 . $$
Equality holds if and only if $\lambda^{ij}=\delta^{ij}$, which means that $\uomega$ is hyperk\"ahler.

The first statement in the Proposition gives a potential obstruction to finding a positive triple with the gives enhanced boundary data. Consider for example the example when  $Y=S^{3}$  and $X$ is the 4-ball. There is a well-known quadratic \lq\lq Chern-Simons'' form $Q_{CS}$ on the closed $2$-forms on $S^{3}$ defined by
$$  Q_{CS}(\mu)= \int_{S^{3}} a \wedge \mu, $$
where $a$ is any $1$-form with $da=\mu$. The necessary condition on our boundary data in this case is that $\mu^{i}$ span a 3-dimensional positive subspace with respect to this form $Q_{CS}$.

\section{Reduction to dimension 3}

We specialise further, mimicking the Gibbons-Hawking construction of hyperk\"ahler $4$-manifolds. Thus we suppose that the $4$-manifold $X$ is the total space of a principal $S^{1}$-bundle over a $3$-manifold $U$, with the action generated by a vector field $\xi$, and consider closed positive triples $\uomega$ which are invariant under the  action. We assume that the action is Hamiltonian for each symplectic structure $\omega^{i}$, so we have Hamiltonian functions $h^{i}:X\rightarrow \bR$ with 
$$  dh^{i}= i_{\xi} \omega^{i}, $$ and these functions are fixed by the  circle action,  so descend to $U$.
More invariantly, writing $\uomega\in  \Omega^{2}(X)\otimes W$ we have $i_{\xi}\uomega\in \Omega^{1}(X)\otimes W$ and this is the derivative of a map $\uh:X\rightarrow W$. The functions $h^{i}$ are then the components of $\uh$ with respect to a co-ordinate system  $W=\bR^{3}$. The definitions imply that $\uh$ induces  a local diffeomorphism from $U$ to $W$, so for local calculations we can suppose that the base $U$ is a domain in $W$ and the functions $h^{i}$ can be identified with standard co-ordinates $x^{i}$ on $W$.  One finds that the general form of such a triple is given by
\begin{equation}  \omega^{i}= \alpha \wedge dx^{i}+ \sum_{{\rm cyclic}} \sigma^{ij} dx^{k} dx^{l}, \end{equation}
where $(jkl)$ run over cyclic permutations,  $\sigma= (\sigma^{ij})$ is a symmetric and positive definite matrix (a function of the co-ordinates $x^{i}$) and $\alpha$ is a connection $1$-form on $X$. The condition that $\sigma$ is symmetric is the same as saying that the connection is the obvious one defined by the metric induced by $\omega^{i}$, with horizontal subspaces the orthogonal complement of $\xi$. 

We will now investigate the reduced $G_{2}$-equations in this context. Write
$F$ for the curvature of the connection, so $F=d\alpha$ and can be regarded as a $2$-form on $U$. We write 
$$F= \sum_{{\rm cyclic}} F^{i} dx^{j} dx^{k}. $$
Now, writing $\partial_{j}$ for partial derivatives,
 $$d\omega^{i} = (F^{i} + \partial_{j} \sigma^{ij}) dx^{1} dx^{2} dx^{3}$$ 
so the condition that we have a closed triple is that 
$$   F^{i}= - \partial_{j}\sigma^{ij}. $$
Now $dF=0$, which is to say
$$    \sum \partial_{i} F^{i}=0, $$
and, at least locally, any closed $2$-form specifies a connection, up to gauge equivalence. So, locally, we can eliminate the connection and curvature and closed triples correspond to matrix-valued functions $\sigma^{ij}$ with
\begin{equation} \sum_{ij} \partial_{i}\partial_{j} \sigma^{ij}= 0. \end{equation}

The volume form  defined by the triple (8) is
$$  \chi  = \det(\sigma)^{1/3} \left(\alpha dx^{1} dx^{2} dx^{3}\right)$$
This implies that
$$ \lambda_{ij}= \sigma_{ij} (\det \sigma)^{1/3}, $$
where $\sigma_{ij}$ denotes the  matrix inverse to $\sigma^{ij}$ as usual. Now the equation (6) is
$$   \sum_{k} \partial_{k}\lambda_{ij} dx^{k}\left( \alpha dx^{j} + \sum_{p,q,r {\rm cyclic}}\sigma^{jp} dx^{q} dx^{r}\right)=0. $$
Expanding this out we get two conditions
\begin{enumerate}\item $\partial_{k}\lambda_{ij}= \partial_{j} \lambda_{ik}$;
\item $\sum_{jk} (\partial_{k}\lambda_{ij}) \sigma^{jk} =0$.
\end{enumerate}
The first condition asserts, at least locally, that $\lambda_{ij}$ is the Hessian of a function, $u$ say:
$$  \lambda_{ij}= \partial_{i}\partial_{j} u. $$
The second condition is implied by the first since 
$$  \partial_{k} \lambda_{ij} \sigma^{jk}= - \partial_{i}\det(\lambda), $$ and the determinant of $\lambda$ is $1$ by construction. To sum up, the 3-dimensional reduction of the $G_{2}$ equations can be written locally as a pair of equations for two functions $u,V$ on a domain in $\bR^{3}$. First, the Monge-Amp\`ere equation \begin{equation}     \det (\partial_{i}\partial_{j} u)=1 \end{equation}
and second 
\begin{equation}       \sum \partial_{i}\partial_{j}(V u^{ij})=0, \end{equation}
where $u^{ij}$ is the inverse of the Hessian $u_{ij}=\partial_{i}\partial_{j} u$.
Given a pair $u,V$ satisfying these equations, we set $\sigma^{ij}= V u^{ij}$ and the discussion above shows that all solutions arise in this way (locally). Notice that, given $u$ the second equation is a linear equation for $V$ and in fact is familiar as  the linearisation  of the Monge-Amp\`ere equation at $u$. Recall that the linearised operator $\Box_{u}$ can be written in three different ways
$$  \Box_{u} f = \sum \partial_{i}\partial_{j}( f u^{ij})= \sum \partial_{i}( \partial_{j}f\  u^{ij})= \sum (\partial_{i}\partial_{j} f) u^{ij}, $$
using the identity $\sum \partial_{i} u^{ij}=0$.

Now we want to set up our boundary value problem in this context. We suppose that $U$ is a $3$-manifold with boundary $\Sigma$ and the circle bundle extends to the boundary, so that $Y=\partial X$  is a circle bundle over $\Sigma$.
 (The extension of the circle bundle over $U$ means that it must be a trivial bundle,  but we do not have a canonical trivialisation.)
 We want to consider triples $\mu^{i}$ of closed $2$-forms on $Y$, as before, invariant under the circle action and such that the action is \lq\lq Hamiltonian'', i.e. there are circle-invariant functions
$  h^{i}_{Y}$ on $Y$ with
   $$  dh^{i}_{Y}= i_{\xi} \mu^{i}, $$
   These functions give a map $\uh_{\Sigma}:\Sigma \rightarrow W=\bR^{3}$ and it follows from the definitions that this is an immersion. Now we encounter a potential obstruction of a differential topological nature to the existence of an invariant  closed positive triple on $X$ with these boundary values: the immersion $\uh_{\Sigma}$ must extend to an immersion of $U$ in $\bR^{3}$. But let us suppose here for simplicity that $\uh_{\Sigma}$ is an embedding of $\Sigma$ as the boundary of a domain in $\bR^{3}$. Then for any extension of the $\rho^{i}$ over $X$, of the kind considered above, the map $\uh$ must be a diffeomorphism from $U$ to this domain. Thus we can simplify our notation by taking $U$ to be a domain in $\bR^{3}$ with smooth boundary $\Sigma$. To avoid complication, we suppose that $U$ is simply connected, so that $\Sigma$ is diffeomorphic to a $2$-sphere.  Thus our PDE problem is to solve the equations (10) and (11) for functions  $u,V$  on $U\subset \bR^{3}$ and the remaining task is to identify the boundary conditions on $\Sigma=\partial U$ defined by a triple $\mu^{i}$. (The assumption that $U$ is simply connected means that the above local analysis of solutions  applies globally on $U$.)

The differential geometric analysis of invariant triples $\umu=(\mu^{i})$ is complicated by the fact that there is no natural connection on the circle bundle $\pi:Y\rightarrow \Sigma$.  But the analysis has a simple conclusion which can be expressed in terms of certain distributions, or currents, which we call {\it layer currents}.  In this analysis it will be important to keep track of the full $SL(3,\bR)$-invariance of the set-up so we work in the $3$-dimensional vector space $W$ with fixed volume element. So  we have an embedding of $\Sigma$ in $W$ as the boundary of a domain $U$ and a triple $\umu$
is a section of $W\otimes \Lambda^{2}T^{*}Y$.

We define a layer current supported on $\Sigma$ to be a linear map from functions on $W$ to $\bR$ of the form 

\begin{equation}  {\cal L}_{\theta_{1}, \theta_{2}, v}(f) = \int_{\Sigma}
 (\nabla_{v} f)\  \theta_{1}+ f \ \theta_{2} . \end{equation}
where $\theta_{1}, \theta_{2}$ are $2$-forms on $\Sigma$ with $\theta_{1}>0$ and $v$ is an outward-pointing normal vector field along $\Sigma$---a section of the tangent bundle of $W$ restricted to $\Sigma$ which is complementary to the tangent bundle of $\Sigma$.
Of course this depends only on the restriction of $f$ to the first formal neighbourhood of $\Sigma$, in particular it is defined for a function $f$ on $U$ which is smooth up to the boundary.

 The point is that the same functional ${\cal L}$ can be defined by different data
$(\theta_{1}, \theta_{2}, v)$. First, it is obvious that for any positive function $g$ on $\Sigma$ the data $(g\theta_{1}, \theta_{2}, g^{-1} v)$ defines the same current. Second,  if $w$ is a tangential vector field on $\Sigma$ we have
\begin{equation}   \int_{\Sigma} (\nabla_{w} f) \theta_{1}= \int_{\Sigma} \tilde{\theta}_{2} \ f \end{equation}
where $\tilde{\theta}_{2} =- d(i_{w} \theta_{1})$. It follows that a given layer current ${\cal L}$ of this mind can be represented using {\it any} normal vector field, for appropriate $\theta_{1}, \theta_{2}$.  Let $\nu_{\Sigma}= TW/T\Sigma$ be the normal bundle of $\Sigma$ in $W$. The fixed volume element on $W$ gives an isomorphism
\begin{equation} \nu_{\Sigma} = \Lambda^{2}T^{*}\Sigma.\end{equation}
Let $[v]$ denote the image of $v$ in $\nu_{\Sigma}$. The product
$$  H_{\cal L} = [v]. \mu_{2} \in \left(\Lambda^{2}T^{*}\Sigma\right)^{\otimes 2} $$
is independent of the choice of data $(v,\mu_{1},\mu_{2})$ used to represent ${\cal L}$; we call $H_{\cal L}$ the {\it primary invariant} of ${\cal L}$. 
For a function $f$ which vanishes on $\Sigma$ the derivative $df$ along $\Sigma$ is defined as a section of $\nu_{\Sigma}^{*}$ and for such functions we have 

\begin{equation}   {\cal L}(f) = \int_{\Sigma} H. df, \end{equation}
where we use the isomorphism $\nu^{*}_{\Sigma}= \Lambda^{2}T^\Sigma$ and the pairing with $H$ yields a $2$-form $H.df$ on $\Sigma$.

Now let $\umu\in W\otimes \Omega^{2}(Y)$ be a closed $S^{1}$-invariant triple on the circle bundle $Y$ over $\Sigma$ such that the inclusion $\Sigma\rightarrow W$ is the Hamiltionian map for the action. Let $y$ be a point of $Y$ and $\epsilon\in W^{*}$ be a co-normal to $\Sigma$ at $x=\pi(y)$, {\it i.e.} an element of $W^{*}$ vanishing on $(T\Sigma)_{x}\subset W$.  Then we have a map
$$    \epsilon: (\Lambda^{2}T^{*}Y)_{y}\otimes W\rightarrow (\Lambda^{2}T^{*}Y)_{y}, $$ and it follows from the definitions that $\epsilon(\umu)$ lies in the image of the pull-back map $\pi^{*}: (\Lambda^{2} T^{*}\Sigma)_{x} \rightarrow (\Lambda^{2}T^{*}Y)_{y}.$ Thus we have a unique element $h\in (\Lambda^{2}T^{*}\Sigma)_{x}$ with $ \pi^{*}(h)= \epsilon(\urho)$. Multiplying $\epsilon$ by a factor $\kappa$ clearly multiplies $h$ by $\kappa$ so, using again the isomorphism (14), we get a well-defined section  $H^{\umu}$ of $(\Lambda^{2}T^{*}\Sigma)^{\otimes 2}$, independent of the choice of $\epsilon$. We call $H^{\umu}$ the primary invariant of the triple $\umu$.

Next choose a normal vector field $v$ along $\Sigma$. At a point $y\in Y$  we transpose $\umu(y)$ to give  a map
 $$ \tilde{\umu}: W^{*}\rightarrow (\Lambda^{2}T^{*}Y)_{y}= TY \otimes \Lambda^{3} T^{*}Y. $$
 The annihiliator of $v(\pi(y))$ is a 2-dimensional subspace of $W^{*}$ and it follows from the definitions that the image of this subspace under $\tilde{\umu}$ defines a a 2-dimensional subspace of $TY$ transverse to the $S^{1}$-orbit. In other words the choice of normal vector field $v$ defines a connection on the $S^{1}$ bundle $\pi:Y\rightarrow \Sigma$: in fact giving a connection is equivalent to giving a complementary bundle to $T\Sigma\subset W$. Let $\Phi$ be the curvature of this connection, a $2$-form on $\Sigma$ and define a current
${\cal L}^{\umu,v}$ by
\begin{equation}  {\cal L}^{\umu,v}(f) =\int_{\Sigma} (H^{\umu}. [v]^{-1})\nabla_{v} f + \Phi f . \end{equation}

Here $H_{\umu}.[v]^{-1}$ is the $2$-form given by the pairing of $[v]^{-1}\in \nu^{*}= (\Lambda^{2}T\Sigma)^{-1}$ with $H^{\umu}\in (\Lambda^{2}T^{*}\Sigma)^{\otimes 2}$. 

\begin{prop}
The layer current ${\cal L}^{\umu,v}$ is independent of the choice of normal vector field $v$ so can be written as ${\cal L}^{\umu}$. Two triples $\umu, \umu'$ are equivalent by $S^{1}$-equivariant diffeomorphisms if and only if ${\cal L}^{\umu}={\cal L}^{\umu'}$.
\end{prop}
 If we change $v$ by multiplication by a positive function  then we do not change the connection and hence we do not change the integral of $\Phi f$. The other term  in the integrand is also unchanged because the scalings of $[v]^{-1}$ and $\nabla_{v}$ cancel. So to prove the first statement it suffices to consider changing $v$ to $v+w$ where $w$ is a tangential vector field on $\Sigma$. Using the formula (13), we have to show that the connection changes by the addition of the $1$-form
$a=i_{w}\theta_{1}$. To see this we work in co-ordinates at a given point on $\Sigma$, taking $v=\partial_{1}$ and the tangent space of $\Sigma$ spanned by $\partial_{2}, \partial_{3}$. Write $\theta_{1}= G dx^{2} dx^{3}$ at the given point. If $\alpha$ is the connection $1$-form on $Y$ defined by $v$ then it follows from the definitions that, over this point,
$$\mu^{1}= G dx^{2} dx^{3}\ , \mu^{2}= \alpha \wedge dx^{2} \ , \ \mu^{3}=\alpha\wedge dx^{3}. $$
If $w= w^{2}\partial_{2}+ w_{3}\partial_{3}$ at this point the annihiliator of $v+w$ in $W^{*}$ is spanned by $dx^{2}- w^{2} dx^{1}, dx^{3}- w^{3} dx^{1}$ and this maps to the $2$-dimensional subspace in $\Lambda^{2}T^{*}Y$ spanned by $$ \alpha \wedge dx^{2}- w^{2} G dx^{2} dx^{3}\ , \ \alpha\wedge dx^{3}- w^{3} G dx^{2} dx^{3}$$
which corresponds to the $2$-dimensional subspace in $TY$ spanned by
$$ \partial_{3}- G w^{2} \xi\ ,\  -\partial_{2}- G w^{3} \xi. $$
This is the anhilliator of the $1$-form $\alpha+a$ where $a= G w^{2} dx^{3}- G w^{3} dx^{2}$ which is the contraction $i_{w} \theta_{1}$ as required.

The second statement of the proposition follows easily from the fact that, since $\Sigma$ is simply connected, a connection is determined  up to gauge equivalence by its curvature.

So far we have considered our structures over the surface $\Sigma\subset W$. Now let $\sigma$ be a matrix-valued function over $U\subset W$, as before, defining a triple $\uomega$ on $X$. Then for any smooth function $f$ on $U$ we have
\begin{equation} \int_{U}\sum \sigma^{ij} \partial_{i}\partial_{j} f - \sum (\partial_{i}\partial_{j}\sigma^{ij}) f  = {\cal L}_{\sigma}(f), \end{equation}
where ${\cal L}_{\sigma}$ is the layer current supported on $\Sigma$ defined by
\begin{equation}
{\cal L}_{\sigma} f = \int_{\Sigma} \sum\sigma^{ij}\partial_{i}f  - (\sum\partial_{i}\sigma^{ij}) f.\end{equation}
(To clarify notation: in (17) we suppress the volume form on $W$ which defines our measure and in (18) the integrand is written as a vector field, which defines a $2$-form on $\Sigma$ by contraction with the 3-dimensional volume form as in (14).)
Then we have: 
\begin{prop}
The boundary value of the triple $\uomega$ corresponding to $\sigma$ is equivalent to  the triple $\umu$ on $Y$ if and only if ${\cal L}_{\sigma}={\cal L}^{\umu}$. \end{prop}

To see this, regard the inverse matrix $\sigma_{ij}$ as a Riemannian metric on $U$. The orthogonal complement with respect to this metric defines a normal vector field $v_{\sigma}$ over $\Sigma$ and hence a connection on $Y\rightarrow \Sigma$.  We know that $\sigma$ defines a connection on the circle bundle $X\rightarrow U$. with curvature given by 
$F^{i}= -\partial_{j} \sigma^{ij}$.  The Proposition amounts to the fact that the restriction of this connection to $Y\rightarrow \Sigma$ is the same as the connection defined by $v_{\sigma}$, which we leave for the reader to check. 

To illustrate the nature of this boundary condition consider an example where $\Sigma$ is locally given by the plane $x^{1}=0$ and take $\partial_{1}$ as normal vector field. Then ${\cal L}^{\umu}$ is locally represented by $2$-forms $$\theta_{1}= G_{1} dx^{2}dx^{3}, \theta_{2}= G_{2} dx^{2} dx^{3}, $$
where $G_{i}$ are functions of $x^{2}, x^{3}$.
That is, for functions $f$ supported in this region
$$  {\cal L}(f)= \int_{x^{1}=0} \left( G_{1} \frac{\partial f }{\partial x^{1}} + G_{2}  f \right)\ dx^{2}dx^{3}. $$
Now if $\sigma$ is defined over $U$ we have, for such functions $f$,
$$  {\cal L}_{\sigma} (f)= \int_{x^{1}=0} \left(\sigma^{11} \frac{\partial f}{\partial x^{1}} + \left(\sigma^{12}\frac{\partial f}{\partial x^{2}}+ \sigma^{13}\frac{\partial f}{\partial x^{3}}\right) - \left(\partial_{i}\sigma^{1i}) f\right)\right)  \   dx^{2}dx^{3}. $$
 Integrating by parts, the sum of the second and third terms is 
 $$ - \int_{x^{1}=0} f (\partial_{1}\sigma^{11}+2\partial_{2}\sigma^{12}+ 2\partial_{3}\sigma^{13}) dx^{2}dx^{3}. $$
Our boundary conditions are
\begin{itemize}
\item  $\sigma^{11}= G_{1}$,
\item $ \partial_{1}\sigma^{11}+2\partial_{2}\sigma^{12}+
2\partial_{3}\sigma^{13}= -G_{2}$. 
\end{itemize}

Notice that if $f$ is an affine-linear function then ${\cal L}_{\sigma}(f)$ vanishes for any $\sigma$ on $U$.  This is connected to the following identities on the boundary:
\begin{itemize}  \item For a circle bundle $Y\rightarrow\Sigma$ with Chern class $d$ and any invariant triple $\umu$ on $Y$, the value of functional ${\cal L}^{\umu}(1)=2\pi d$ (Here $1$ denotes the constant function). .
\item Suppose $d=0$, so $Y$ is diffeomorphic to $S^{1}\times \Sigma=S^{1}\times S^{2}$ and there is a lift
$[\Sigma]\in H_{2}(Y)$. Then for any invariant triple $\urho$ 
$$   {\cal L}^{\umu}(x^{i})= \int_{[\Sigma]} \mu^{i}. $$
\end{itemize}
 
 Again, we leave the proofs as exercises for the interested reader.

Putting all this together, we can formulate the dimensionally-reduced version of our general  boundary value problem as follows. The functional (9) clearly reduces to the functional 
\begin{equation}  {\rm Vol}(\sigma) = \int_{U} (\det\sigma) ^{1/3}. \end{equation}

{\bf Variational Problem I}

{\it Given a (simply connected) domain $U\subset \bR^{3}$ with smooth boundary $\Sigma$ and a layer current  ${\cal L}$ on $\Sigma$, find the critical points of the volume functional (19) over all $\sigma =(\sigma^{ij})$ on $U$ satisfying}
\begin{itemize}
\item (A) \ $\sum \partial_{i}\partial_{j}\sigma^{ij}=0$,
\item (B)\ ${\cal L}_{\sigma} = {\cal L}$.
\end{itemize}

Our first question is now whether the set ${\cal S}_{{\cal L}}$ of matrix-valued functions $\sigma$ satisfying (A), (B) above is non-empty. The integral formula gives an immediate  constraint on the boundary data: if ${\cal S}_{{\cal L}}$ is not empty then ${\cal L} f\geq 0$ for all convex functions $f$ on $U$ (with equality if and only if $f$ is affine linear). We also have
\begin{prop}
If a solution to the Variational Problem I exists it is an absolute maximum of the  volume functional  on ${\cal S}_{{\cal L}}$.
\end{prop}

This  follows immediately  from the facts that  both conditions (A),(B) are linear in $\sigma$ and the function $(\det\sigma)^{1/3}$ is concave.

Let $f$ be a convex function on $U$ satisfying the Monge-Amp\`ere equation
${\rm det} (f_{ij})=1$. Then for any $\sigma$ we have, pointwise on $U$, \begin{equation}   \det(\sigma)^{1/3} = \left(\det \sigma_{ij} \det(f_{ij})\right)^{1/3}\leq \frac{1}{3} \sum \sigma^{ij}f_{ ij}. \end{equation}
So if $\sigma$ satisfies the conditions (1),(2) of Variational Problem I we have, integrating over $U$ and using the definition of ${\cal L}_{\sigma}$, 
\begin{equation} {\rm Vol}(\sigma) \leq\frac{1}{3} \int_{U} \sum \sigma^{ij} f_{ij}= \frac{1}{3} {\cal L}_{\sigma} f. \end{equation} 

(Our previous bound, in Proposition 2, arises by taking quadratic functions $f$. )
These bounds furnished by solutions of the Monge-Amp\`ere equation, lead to a dual formulation of the variational problem, which incorporates the boundary conditions in a simple way. Write $MA(U)$ for the set of convex solutions of the Monge-Amp\`ere equation on $U$, smooth up to the boundary.

\

{\bf Variational problem II}

{\it Given a (simply connected) domain $U\subset \bR^{3}$ with smooth boundary
$\Sigma$ and a layer current ${\cal L}$ on $\Sigma$, minimise ${\cal L}(f)$ over all $f\in MA(U)$.}

\begin{prop}

The variational problems I,II are equivalent in the sense that for $u\in MA(U)$ we can find a positive function $V$ such that $\sigma^{ij}= V u^{ij}$ is a solution of variational problem I if and only $u$ is a solution of variational problem II.
\end{prop}

 In one direction, equality holds in (21) if and only if $u_{ij}$ is a multiple of the inverse of $\sigma^{ij}$. We know that a solution to the variational problem I has the form
$\sigma_{ij}= V u^{ij}$ where $u$ satisfies the Monge-Amp\'ere equation, so taking $f=u$ equality holds in (21),  and  it follows that $u$ minimises ${\cal L}(f)$ over $MA(U)$. 

In the other direction, suppose that $u\in MA(U)$ is an extremum of the functional ${\cal L}$ (it will follow from the discussion below that $u$ is in fact a minimum and is unique up to the addition of an affine-linear function). The Euler-Lagrange equation is  ${\cal L}(W)=0$ for all solutions $W$ of the linearised equation $\Box_{u}W$ over $U$. We can solve the Dirichlet problem for this linearised equation, to find $V$ such that $\Box_{u}V=0$ and so that if $\sigma^{ij}= Vu^{ij}$ the  primary invariant of ${\cal L}_{\sigma}$ is equal to that of ${\cal L}$.  If we chose a co-normal $\nu^{*}$ to $\Sigma$ this is just saying that   $\sum Vu^{ij} \nu^{*}_{i}\nu^{*}_{j}$ is a prescribed function on $\Sigma$, which for fixed $u$ is just prescribing $V$ on $\Sigma$. Then it follows from the previous discussion that ${\cal L}_{\sigma}(W)=0$ for all solutions $W$ of the linearised equation. Since ${\cal L}_{\sigma}$ and ${\cal L}$ have the same primary invariant so the difference can be written as
$$  ({\cal L}_{\sigma}-{\cal L})(f)= \int_{\Sigma} \Theta f, $$
for a $2$-form $\Theta$ on $\Sigma$. For any function $f$ on $\Sigma$ we can solve the Dirichlet problem for $\Bix_{u}$ with boundary value $f$ and so $$  \int_{\Sigma} \Omega f =0$$ for all $f$. This implies that  $\Omega=0$ so ${\cal L}_{\sigma}={\cal L}$ and we have solved the variational problem I. 

\

 Modifying our problem, we can obtain a decisive existence result. Rather than fixing the full boundary data $\umu$  we just fix the primary invariant $H^{\umu}$. Given a positive $H\in \Gamma(\Lambda^{2}T^{*}\Sigma)^{\otimes 2}$ we write  ${\cal C}_{H}$ for the set of $\sigma$ over $U$ satisfying $\sigma^{ij}_{,ij}=0$ and with the primary invariant of ${\cal L}_{\sigma}$ equal to $H$. As above, in terms of a co-normal $\nu^{*}$ this amounts to prescribing
$\nu^{*}_{i}\nu^{*}_{j}\sigma^{ij}$ on the boundary. 
\begin{prop}
If $U$ is strictly convex  there is a unique critical point of the volume functional  on ${\cal C}_{H}$ and this is an absolute maximum.
\end{prop}
The uniqueness and the fact that a  critical point is an absolute maximum follows from concavity,  just as before. For the existence,  we first solve
(invoking \cite{kn:CNS}) the Dirichlet problem for the Monge-Amp\`ere equation to get a function $u\in MA(U)$ with $u=0$ on $\Sigma$.
\ Now solve the Dirichlet problem for the linearised equation to find a function $V$ with $ \Box_{u}=0$ in $U$ and such that $V u^{ij} \nu^{*}_{i}\nu^{*}_{j}$ is the prescribed function on the boundary and write $\sigma^{ij}=V u^{ij}$. We claim that  this $\sigma$ is a critical point of the volume functional on ${\cal C}_{H}$. Let $\tau^{ij}$ be an infinitesimal variation within ${\cal C}_{H}$. In other words, $\sum \partial_{i}\partial_{j}\tau^{ij}=0$ in $U$ and on the boundary  $\sum \tau^{ij} \nu^{*}_{i}\nu^{*}_{j}=0$.  Then
the variation in the volume functional is
\begin{equation} 3\delta {\rm Vol}= \int_{U} \sum u_{ij} \tau^{ij} = \int_{\Sigma}
 \sum \tau^{ij}_{j} u - \sum \tau^{ij} \partial_{j}u. \end{equation}
 
The first term on the right hand side of (22) vanishes since $u$ vanishes on $\Sigma$. In the  second term, the derivative of $u$ along $\Sigma$ vanishes,  so there is only a contribution from the normal derivative of $u$ and the integrand is a multiple of $\sum \tau^{ij} \nu^{*}_{i}\nu^{*}_{j}$,  so this also vanishes.

\section{Further remarks}
\subsection{Singularities}

It seems unlikely that the variational problems I,II always have solutions, even given the constraints we have found. To see this we consider the well-known singular solutions of the Monge-Amp\`ere equation, going back to Pogerolov. With co-ordinates $x^{1}, x^{2}, x^{3}$ set $r= \sqrt{(x^{1})^{2}+ (x^{2})^{2}}$ and consider functions $u$ of the form $u= f(x^{1}) r^{4/3}$. Then one finds that 
$$  \det(u_{ij})= \frac{64}{27} f \left( \frac{f f''}{3}- (f')^{2}\right), $$
so we can find smooth functions $f$ on an interval, say $(-\epsilon, \epsilon)$ with $f''>0$ and such that $u$ satisfies the Monge-Amp\`ere equation. Fix such a function $f$ and let $\eta$ be the vector field
$$ \eta= 2x^{1} \frac{\partial}{\partial x^{1}} - x^{2} \frac{\partial}{\partial x^{2}}- x^{3} \frac{\partial}{\partial x^{3}}. $$
This vector field generates volume preserving transformations, so $V=\nabla_{\eta} u$ satisfies the linearised equation $\Bix_{u} V=0$ and if we define $\sigma^{ij}= V u^{ij}$ we get a singular solution of our reduced $G_{2}$ equations (provided that $V>0$).  Suppose that, near the origin, $\Sigma$ is given by the co-ordinate plane $x^{1}=0$ and let our boundary data be given locally by the layer current
$$   \int G_{1} \partial{f}{\partial x^{1}} _{1} +  G_{2} f  dx^{2} dx^{3}, $$
 as above, for smooth functions $G_{i}(x^{2}, x^{3})$. As we saw above, the boundary condition is given by
$$   \sigma^{11} = G_{1}  \  ,  \   \partial_{1} \sigma^{11} +2\partial_{2} \sigma^{12}+2\partial_{3} \sigma^{13} = - G_{2}. $$
 
 One can compute that
$$  V= (2x^{1} f' - \frac{4}{3} f) r^{4/3} $$
$$   \sigma^{11}= \frac{16 f^{2}}{9} ( 2x^{1} f'- \frac{4}{3}f) $$
$$   \sum_{i} \partial_{i} \sigma^{1i}= \frac{16 f^{2}}{9} (2x^{1} f' - \frac{4}{3} f)'- \frac{16}{3} f f'.  $$
Thus $\sigma^{11}, \sum_{i} \partial_{i}\sigma^{1i}$ are smooth functions of $x^{1}$ so  we get a singular solution of our boundary value problem with smooth boundary data $G_{1}, G_{2}$. 

\subsection{Connection with the Apostolov-Salamon construction}

In the discussion above we have passed from 7 dimensions to 3 dimensions by first imposing translational symmetry in 3 variables to get down to 4 dimensions and then imposing a circle action to pass from 4 to 3. We can achieve the same end by imposing the circle action first, to get a reduction to 6 dimensions, and then studying translation invariant solutions. The material in this subsection was explained to the author by Lorenzo Foscolo. 

$G_{2}$ structures on a $7$-manifold $M$ invariant under a free circle action have been studied by Apostolov and Salamon \cite{kn:AS} and others. The quotient space $N$ has an induced $SU(3)$ structure, that is to say a $2$-form $\omega$ and a complex $3$-form $\Omega$ equivalent at each point to the standard structures on $\bC^{3}$ (with complex $3$-form $dz_{1}dz_{2} dz_{3}$). The $G_{2}$ structure on $M$ can be written as
\begin{equation}
\phi = \alpha \wedge \omega + V^{3/4} {\rm Re}\  \Omega \end{equation}
where we identify forms on $N$ with their lifts to $M$ and $V$ is a smooth positive function on $N$. (In fact $V^{-1/2}$ is the length of generator of the circle action in the metric $g_{\phi}$.) The $1$-form $\alpha$ is a  connection form on the circle bundle $M\rightarrow N$. Now one finds that
$$  *_{g_{\phi}}\phi= -V^{1/4} \alpha\wedge {\rm Im}\ \Omega + \frac{V}{2}\omega^{2}. $$
Thus if $F=d\alpha$ is the curvature of the connection the conditions to be satisfied for a torsion free $G_{2}$-structure are:
\begin{equation} d\omega=0\  , \ F\wedge \omega+ d(V^{3/4} {\rm Re}\ \Omega)=0 \end{equation}
and
\begin{equation} d(V^{1/4} {\rm Re} \Omega)=0 \ ,\ dV\wedge \omega^{2}= 2 V^{1/4} F\wedge {\rm Im}\ \Omega. \end{equation}

Now let $W$ be a 3-dimensional real vector space as before, and set $N=U\times W^{*}$ where $U$ is an open set in $W$. Take standard co-ordinates $x^{i}$ on $W$ and $\theta_{i}$ on $W^{*}$ so there is a standard symplectic form
$$\omega = \sum dx^{i}\wedge d\theta_{i}. $$
Let $u$ be a convex function on $U$ and define complex $1$-forms
$$  \epsilon_{a} =d\theta_{a}+ i \sum u_{ab} dx^{b}. $$
It is well-known that these define a complex structure compatible with $\omega$ and  with holomorphic $3$-form $ \epsilon_{1}\epsilon_{2}\epsilon_{3}$. If $u$ satisfies the Monge-Amp\`ere equation $\det u_{ij} = 1$ then this is  a Calabi-Yau structure. If $V$ is a positive function on $U$ we modify this by taking
 $$  \epsilon'_{a}= V^{-1/4} d\theta_{a}+ i V^{1/4} \sum u_{ab}dx^{b}, $$
 and $$\Omega=\epsilon'_{1}\epsilon'_{2}\epsilon'_{3}.$$
 This complex $3$-form is also algebraically compatible with $\omega$.  One checks that if $u$ satisfies the Monge-Amp\`ere equation and $F$ is the $2$-form
  $$   F= -\sum_{ j,k,l \ {\rm cyclic}}  (\partial_{i}) V u^{ij} dx^{k}dx^{l}$$ then $V, F,\Omega, \omega$ satisfy the equations (24),(25). The condition that $F$ is closed, so arises as the curvature of a connection,  is the equation
$  \Box_{u} V=0 $.

\subsection{A general class of equations and LeBrun's construction}

 Our variational problem I, and its dual formulation II have natural extensions. We can clearly replace $\bR^{3}$ by $\bR^{n}$,  but more interestingly we can consider a class of different functionals. Let $W$ be an $n$-dimensional real vector space and write $s^{2}_{+}(W)\subset s^{2}(W)$ for the cone of positive definite quadratic forms on $W$. Let $\nu$ be a smooth positive concave function on $s^{2}_{+}(W)$ which is homogeneous of degree $1$, so
$\nu(k \sigma)= k \nu(\sigma)$. This means that $\nu$ is not strictly concave but we suppose that the kernel of the Hessian of $\nu$ has dimension $1$ everywhere which means that $\log \nu$ is strictly concave. The basic example, which we have discussed in the case $n=3$,  is to take $\nu(\sigma)= (\det \sigma)^{1/n}$. In the general case we consider  functions $\sigma$ on a simply-connected domain $U\subset W$ taking values in $s^{2}_{+}(W)$, which we express in terms of a basis as $\sigma=(\sigma^{ij})$. We consider the functional
$$ I=\int  \nu(\sigma) $$
subject to the constraint  $\sum \partial_{i}\partial_{j}\sigma^{ij}=0$. Initially we consider   variations of this integral with respect to compactly-supported variations of $\sigma$, and later we consider boundary conditions. 

 We regard
$D\nu$ as a map from $s^{2}_{+}(W)$ to $s^{2}(W^{*})$. The homogeneity of $\nu$ implies that this map is constant on rays in $s^{2}_{+}(W)$, so it cannot map  onto an open set in  $s^{2}(W^{*})$. Let $\hat{L}$ be the Legendre transform of $L=\log \nu$.  So $\hat{L}$ is a concave function on some open subset $\Omega$ of $s^{2}(W^{*})$ and $ \hat{N}=\hat{L}^{-1}(1)$ is a smooth hypersurface in $\Omega$. The homogeneity of $\nu$ implies that the image of $D\nu$ is exactly $V$ and for $v\in \hat{N}$ the pre-image $(D\nu)^{-1}(v)$ is a single ray in $s^{2}_{+}(W)$. Set $N=\nu^{-1}(1)\subset s^{2}_{+}(W)$. Then the restriction of $D\nu$ gives a diffeomorphism from $N$ to $\hat{N}$ and we write $\psi:\hat{N}\rightarrow N$ for the inverse.   
\begin{prop}
The integral $I$ is stationary, with respect to compactly-supported variations in $\sigma$ satisfying the constraint $\sum \partial_{i}\partial_{j}\sigma^{ij}=0$, if and only if there is a function $u$ on $U$ such that 
$$  \partial_{i}\partial_{j} u = (D\nu)(\sigma). $$
\end{prop}
In one direction, if there is such a function $u$ and if $\tau$ is a compactly supported variation of $\sigma$ with $\partial_{i}\partial_{j} \tau^{ij} =0$ then
$$ \delta I = \int \langle D\nu(\sigma), \tau\rangle =\int (\partial_{i}\partial_{j} u) \tau^{ij} = \int u \partial_{i}\partial_{j} \tau^{ij} = 0. $$
The other direction follows easily from the fact that we can generate solutions of $\partial_{i}\partial_{j} \tau^{ij}=0$ from an arbitrary tensor $h^{iaj}$ which is skew symmetric in $i,a$   via the formula
$$   \tau^{ij}= \partial_{a} h^{iaj} + \partial_{a} h^{jai}. $$

Now the general local solution of our variational problem is obtained as follows. First solve the equation for a function $u$ that $u_{ij}=\partial_{i}\partial_{j} u $ lies in $\hat{N}$; that is
\begin{equation} \hat{L}(u_{ij}) =1. \end{equation}
Now set $\lambda=\left(\lambda^{ij}\right)= \psi(u_{ij})$. We know that $\sigma = V \psi(u_{ij})$ for some positive function $V$ and the remaining equation to solve is the linear equation for $V$
\begin{equation}    \partial_{i}\partial_{j} ( \lambda^{ij} V )=0. \end{equation}
We compare this with the linearisation of the nonlinear equation (26) at $u$. The derivative of $\hat{L}$ at a point $v\in N$ is given by $\psi(v)$, so the linearised equation is 
\begin{equation} \sum \lambda^{ij} \partial_{i}\partial_{j} V=0. \end{equation}
In general this is not the same as the equation (27), but the two equations are adjoint in that the formal adjoint $\Diamond_{u}^{*}$ of the operator
$\Diamond_{u}(V)= \sum \lambda^{ij}\partial_{i}\partial_{j}V $ is $\Diamond^{*}(V)= \sum \partial_{i}\partial_{j}(\lambda^{ij} V)$. 

In the case when $\nu(\sigma)= (\det \sigma)^{1/n}$ the derivative of $\log \nu$ is the map $\sigma\mapsto n^{-1}\sigma^{-1}$ and we  recover the previous set-up. In this special case the linearised equation is self-adjoint, i.e $\Diamond^{*}_{u}=\Diamond_{u}$. 

We can now introduce a boundary value problem on  a domain $U\subset W=\bR^{n}$ with smooth boundary $\Sigma$ and with a given layer current${\cal L}$ supported on $\Sigma$, extending the definitions from $n=3$ in the obvious way. If $u$ solves the nonlinear equation $\hat{L}(u_{ij})=1$ on $U$ then for any $\sigma$ satisfying the boundary conditions and $\sum \partial_{i}\partial_{j}\sigma^{ij}=0$ we have an inequality 
$$  \int_{U} \nu(\sigma) \leq {\cal L}(u), $$ and we get a dual variational problem as before. 

 Claude LeBrun pointed out to the the author that there are some striking similarities between the variant of the Gibbons-Hawking construction studied in the previous section and another variant introduced by him in \cite{kn:Leb}, constructing K\"ahler surfaces of zero scalar curvature. While this does not exactly fit into the general framework above, we will outline how it can be treated in a similar fashion. 

We consider a triple of forms $\omega^{i}$ with $\omega^{1}$ a K\"ahler form and $\omega^{2}, \omega^{3}$ the real and imaginary parts of a holomorphic $2$-form. This means that we restrict attention to matrix-valued functions $\sigma$ which are diagonal, with $\sigma^{11}=a$ and $\sigma^{22}=\sigma^{33}=b$ for positive functions $a,b$. The condition $\partial_{i}\partial_{j}\sigma^{ij}=0$ is then
\begin{equation} a_{11} + b_{22} + b_{33}=0, \end{equation}
(writing $a_{11}=\partial_{1}\partial_{1} a$ etc.).
The variational formulation,  generating the zero scalar curvature equation, comes from the \lq\lq Mabuchi functional'', which in this situation is given by
\begin{equation} I = \int a (\log(a/b)-1), \end{equation}
and the function $\nu(a,b)= a (\log(a/b)-1)$ is homogeneous of degree $1$. The condition that $I$ is stationary with respect to compactly supported variations satisfying the constraint (29) is that there is a function $u$ with $$u_{11}= \log(a/b)\ \ \ u_{22} + u_{33} = -(a/b). $$
In other words, $u$ satisfies the nonlinear equation
\begin{equation}  e^{u_{11}} + u_{22} + u_{33} =0. \end{equation}

Given such a function $u$ we set $b=V$ and $a= e^{u_{11}} V$ and the equation (29) is the linear equation for $V$:
\begin{equation} (e^{u_{11}} V)_{11} + V_{22} + V_{33}=0. \end{equation}
This is again the adjoint of the linearisation of the nonlinear equation (31). To relate this to LeBrun's set-up we put $U= u_{11}$ so
$$ \left( e^{U}\right)_{11}+ U_{22}+U_{33}=0; $$
 $$ (e^{U}V)_{11}+ V_{22}+ V_{33}=0, $$  and these are the equations, for functions $U,V$,  obtained by LeBrun.



\end{document}